\def\version{\today}

\documentclass[reqno,twoside,11pt]{amsart}
\usepackage{cite}
\usepackage{amsmath}
\usepackage{amsfonts}
\usepackage{amssymb}
\usepackage{epsfig}
\usepackage{verbatim}
\usepackage{color}

\usepackage{ulem}

\usepackage{mathpazo}
\usepackage{mathabx}
\usepackage{mathrsfs}

\DeclareFontFamily{OT1}{cmss}{} \DeclareFontShape{OT1}{cmss}{m}{n} {<5> <6> <7> <8> <9> <10> <11> <12> <13> <14.4> cmss10}{}
\DeclareMathAlphabet{\cmss}{OT1}{cmss}{m}{n}

\DeclareFontFamily{OT1}{fraktura}{}
\DeclareFontShape{OT1}{fraktura}{m}{n} {<5> <6> <7> <8> <9> <10> <11> <12> <13> <14.4> [1.1] eufm10}{}
\DeclareMathAlphabet{\fraktura}{OT1}{fraktura}{m}{n}

\newtheoremstyle{thm}{1.5ex}{1.5ex}{\itshape\rmfamily}{} {\bfseries\rmfamily}{}{2ex}{}

\newtheoremstyle{def}{1.5ex}{1.5ex}{\rmfamily\sl}{} {\bfseries\rmfamily}{}{2ex}{}

\newtheoremstyle{rem}{1.3ex}{1.3ex}{\rmfamily}{} {\bfseries\rmfamily}{}{2ex}{}

\newtheoremstyle{ass}{1.5ex}{1.5ex}{\rmfamily\sl}{} {\bfseries\rmfamily}{}{2ex}{}

\newenvironment{proofsect}[1] {\vskip0.1cm\noindent{\rmfamily\itshape#1.}}{\qed\vspace{0.15cm}}

\theoremstyle{thm}
\newtheorem{theorem}{Theorem}[section]
\newtheorem{lemma}[theorem]{Lemma}
\newtheorem{proposition}[theorem]{Proposition}

\newtheorem*{Main Theorem}{Main Theorem.}

\newtheorem{assumption}[theorem]{Assumption}

\newtheorem{definition}[theorem]{Definition}

\theoremstyle{rem}
\newtheorem{remark}[theorem]{{Remark}}

\numberwithin{equation}{section}

\renewcommand{\section}{\secdef\sct\sect}
\newcommand{\sct}[2][default]{\refstepcounter{section}
\addcontentsline{toc}{section}
{{\tocsection {}{\thesection}{\!\!\!\!#1\dotfill}}{}}
\vspace{0.7cm}
\centerline{ 
\scshape\arabic{section}.\ #1} \nopagebreak \vspace{0.2cm}}
\newcommand{\sect}[1]{
\vspace{0.4cm} \centerline{\large\scshape\rmfamily #1}
\vspace{0.2cm}}

\renewcommand{\subsection}{\secdef\subsct\sbsect}
\newcommand{\subsct}[2][default]{\refstepcounter{subsection}
\addcontentsline{toc}{subsection}
{{\tocsection{\!\!}{\hspace{1.2em}\thesubsection}{\!\!\!\!#1\dotfill}}{}}
\nopagebreak\vspace{0.45\baselineskip} {\flushleft\bf
\thesection.\arabic{subsection}~\bf #1.~}
\\*[3mm]\noindent
\nopagebreak}
\newcommand{\sbsect}[1]{
\vspace{0.1cm}\noindent
\textbf{#1.~}\vspace{0.1cm}}

\renewcommand{\subsubsection}{%
\secdef \subsubsect\sbsbsect}
\newcommand{\subsubsect}[2][default]{%
\refstepcounter{subsubsection} 
\addcontentsline{toc}{subsubsection}{{\tocsection{\!\!}
{\hspace{3.05em}\thesubsubsection}{\!\!\!\!#1\dotfill}}{}}
\nopagebreak
\vspace{0.15\baselineskip} \nopagebreak {\flushleft\rmfamily
\itshape\arabic{section}.\arabic{subsection}.\arabic{subsubsection}
\ \rmfamily #1\/.}\ }
\newcommand{\sbsbsect}[1]{\vspace{0.1cm}\noindent
\rmfamily \itshape
\arabic{section}.\arabic{subsection}.\arabic{subsubsection} \
\sffamily #1\/.\ }

\renewcommand{\caption}[1]{%
\vglue0.5cm
\refstepcounter{figure}
\begin{center}
\begin{minipage}[c]{0.8\textwidth}\small {\sc Fig.~\thefigure\ }#1\end{minipage}
\end{center}
}

\newcommand{\textd}{\text{\rm d}\mkern0.5mu}

\newcommand{\texte}{\text{\rm  e}\mkern0.7mu}

\newcommand{\FF}{\mathcal F}

\newcommand{\HH}{\mathcal H}

\newcommand{\LL}{\mathcal L}

\newcommand{\E}{\mathbb E}

\newcommand{\BbbP}{\mathbb P}
\newcommand{\Q}{\mathbb Q}
\newcommand{\R}{\mathbb R}

\newcommand{\Z}{\mathbb Z}

\def\myffrac#1#2 in #3{\raise 2.6pt\hbox{$#3 #1$}\mkern-1.5mu\raise 0.8pt\hbox{$#3/$}\mkern-1.1mu\lower 1.5pt\hbox{$#3 #2$}}

\newcommand{\wt}{\widetilde}

\newcommand{\myemph}[1]{\textit{#1}}

\newcommand\ON{{\text{\,\tiny\rm ON}}}
\newcommand\OFF{{\text{\,\tiny\rm OFF}}} 

\begin{document}

\title[Degenerate 1D random walk\hfill]{An invariance principle for one-dimensional random walks\\in degenerate dynamical random environments}
\author[\hfill M.~Biskup and M.~Pan]
{Marek~Biskup and Minghao Pan}
\thanks{\hglue-4.5mm\fontsize{9.6}{9.6}\selectfont\copyright\,\textrm{2022}\ \ \textrm{M.~Biskup and M.~Pan.
Reproduction, by any means, of the entire
article for non-commercial purposes is permitted without charge.\vspace{2mm}}}
\maketitle

\vspace{-5mm}
\centerline{
\textit{Department of Mathematics, UCLA, Los Angeles, California, USA}}
\smallskip

\smallskip 
\centerline{\version}

\vskip0.5cm
\begin{quote}
\footnotesize \textbf{Abstract:}
We study random walks on the integers driven by a sample of time-dependent nearest-neighbor conductances that are bounded but are permitted to vanish over time intervals of positive Lebesgue-length.  Assuming only ergodicity of the conductance law under space-time shifts and a moment assumption on the time to accumulate a unit conductance over a given edge, we prove that the walk scales, under a diffusive scaling of space and time, to a non-degenerate  Brownian motion for a.e.\ realization of the environment. The conclusion particularly applies to random  walks  on one-dimensional dynamical percolation subject to  fairly  general stationary edge-flip dynamics.
\end{quote}

\section{Definitions and main results}
\noindent
 This note is concerned with large-scale behavior of a particular class of one-dimensional  nearest-neighbor random walks in dynamical random
environments. Each of our random walks is technically a continuous-time Markov chain on~$\Z$ with time-varying generator $L_t$ 
at time~$t$ that acts on $f\colon\Z\rightarrow \R$ as
\begin{equation}
\label{E:1.2}
L_{t}f(x):=\sum_{z=\pm 1}a_{t}(x,x+z)\left[\, f(x+z)-f(x)\right].
\end{equation}
Here the coefficients $a_{t}(x,x+z)$ are non-negative numbers with the intuitive meaning of the jump rate from~$x$ to~$x+z$ at time~$t$. The key restriction we impose is that this jump rate is symmetric,
\begin{equation}
a_{t}(x,x+z)=a_{t}(x+z,x),\quad x\in\Z,\,z=\pm 1,
\end{equation}
and so $a_t(e)$ is just a function of the unordered edge~$e$.
No jump across edge~$e$ can occur at time~$t$ when $a_t(e)$ vanishes.

In order to construct the  Markov  chain precisely we need to make some regularity assumptions on the environment. 
 Writing~$E(\Z)$ for the set of unordered edges of~$\Z$, let $\Omega :=[0,\infty )^{\R\times E(\mathbb{Z)}}$ denote  the set of all environments  and   $\mathcal{F}:=\bigotimes _{\R\times E(\mathbb{Z)}}\mathcal{B([}0,\infty \mathcal{))}
$  for  the product $\sigma $-algebra on~$\Omega$. 
For each~$t\in\R$ and $x\in\Z$ let $\tau _{t,x}\colon\Omega \rightarrow \Omega$
be the  canonical  space-time shift acting  on $a\in\Omega$  as
\begin{equation}
(\tau_{t,x}a)_{s}(y,y+z)=a_{t+s}(y+x,y+x+z),\quad s\in \R%
,\,  x  \in \Z,\, z=\pm1\text{.}
\end{equation}
We will assume throughout that a probability measure~$\BbbP$ on $(\Omega,\FF)$ is given, with expectation denoted as~$\E$, such that the following holds:

\begin{assumption}
\label{ass-1}
For each edge $ e \in E(\Z)$, the map $t \rightarrow a_{t}(e)$
 is Borel measurable and locally Lebesgue integrable. The law $\mathbb{P}$ is
invariant and ergodic with respect to the family of space-time shifts $\left\{ \tau _{t,x}\colon t\in \R%
,\,x\in \Z\right\}$.
\end{assumption}

Under Assumption~\ref{ass-1}, a $\Z$-valued Markov chain with generator \eqref{E:1.2} can be constructed for all environments in a measurable set~$\Omega_0$ of full $\BbbP$-measure. (See~\cite{Bi1} for an outline of that construction with non-explosivity being its main concern.) Let $X=\left\{ X_{t}\colon t\geq 0\right\} $ denote the c\`adl\`ag trajectory of the chain and write $P_{a }^{x}$ to denote the law of~$X$ in environment $a \in \Omega_0$ subject to the initial condition $P_{a }^{x}(X_{0}=x)=1$. 
The aim of the present note is to give sufficient conditions under which the walk behaves ``usually'' at large space-time scales. We formalize this as:

\begin{definition}
We say that a Quenched Invariance Principle holds if there exists a constant $%
\sigma^2\in(0,\infty) $ such that for any $t_0>0$ and $\mathbb{P}$-a.e.\ environment $a $%
, the law of 
\begin{equation}
t\mapsto\frac{1}{\sqrt{n}}X_{nt},\quad\text{ \ }0\leq t\leq t_0,
\end{equation}%
induced by $P_{a }^{0}$ on the Skorohod space $D[0,t_0]$ of c\`{a}dl\`{a}%
g paths converges, as $n\rightarrow \infty $, to the law of Brownian motion $%
\left\{ B_{t}\colon t\in[0,t_0]\right\} $ with $EB_{t}=0$ and $EB_{t}^{2}=\sigma ^{2}t$.
\end{definition}

We note that, for one-dimensional walks subject to Assumption~\ref{ass-1}, a Quenched Invariance Principle was proved earlier by Deuschel and Slowik~\cite{DS16} assuming finiteness of $p$-th positive and $q$-th negative moments of $a_t(e)$ subject to $p,q\ge1$ and $\frac1{p-1}(1+\frac1q)<1$. The latter inequality stems from the method of proof, which is based on elliptic regularity techniques. In~\cite{Bi1}, the first author discovered a different proof that works solely under the first-moment conditions 
\begin{equation}
\label{E:1.4}
\mathbb{E}[a_{t}(e)]<\infty \quad\text{ and }\quad\mathbb{E}[a_{t}(e)^{-1}]<\infty .
\end{equation}
These were also shown to be necessary for the result to hold in general. 

Unfortunately,  under  Assumption~\ref{ass-1},  the negative moment condition in \eqref{E:1.4}  makes it impossible for $t\mapsto a_t(e)$ to vanish on a set of positive Lebesgue measure. This excludes natural examples of prime interest. We mend this partially in:

\begin{theorem}
\label{main}In addition to Assumption~\ref{ass-1}, suppose that
\begin{enumerate}
\item[(1)] $a_{t}(e)\in \lbrack 0,1]$ for all $t\in\R$ and $e\in\E(\Z)$,
\item[(2)] the quantity
\begin{equation}
\label{E:1.4a}
T:=\inf \left\{ t\geq 0:\int_{0}^{t}a_{t}(0,1)\,\textd t\geq 1\right\}
\end{equation}
obeys
\begin{equation}
\label{E:1.6}
\exists\varepsilon>0\colon\quad \mathbb{E}(T^{3+\varepsilon })<\infty .
\end{equation}
\end{enumerate}
Then a Quenched Invariance Principle holds.
\end{theorem}


An important  family  of examples covered by Theorem~\ref{main}, but not the  conclusions  of~\cite{DS16,Bi1}, are random walks on \myemph{dynamical percolation}. Here $a_t(e)$ takes values in~$\{0,1\}$, with value~$1$ representing the edge being ``ON'' and~$0$  for  the edge being ``OFF.'' The processes $\{t\mapsto a_t(e)\}_{e\in\E(\Z)}$ are i.i.d.\ copies of a given stationary process on~$\{0,1\}$  which we assume to have  c\`adl\`ag (and thus piecewise-constant) sample paths and take both values a.s.  To make a connection to percolation we note that, at each given time~$t\in\R$, the configuration of the ``ON'' edges is Bernoulli with probability $p:=\E a_0(0,1)$.

While the nature of the individual edge dynamics can be quite arbitrary,  the assumptions permit a representation via  a sequence of pairs of strictly positive random variables
\begin{equation}
\label{E:1.8a}
\bigl\{(T_i^\OFF,T_i^\ON)\bigr\}_{i\in\Z},
\end{equation}
to be called ``OFF'' and~``ON''-times, that stand for the lengths of successive time intervals on which $t\mapsto a_t(e)$ equals~$0$ and~$1$, respectively.  Explicitly, writing $\{\tau_i\}_{i\in\Z}$ for the successive times when~$t\mapsto a_t(e)$ switches from~$1$ to~$0$ and $\{\tau_i'\}_{i\in\Z}$ for the times it switches from~$0$ to~$1$, indexed so that $\tau_i<\tau_i'<\tau_{i+1}$ for each~$i\in\Z$ and $\tau_0<0<\tau_1$ a.s., these are defined as $T_i^\OFF:=\tau_i'-\tau_i$ and $T_i^\ON:=\tau_{i+1}-\tau_i'$.

The sequence \eqref{E:1.8a} in turn determines the trajectory~$t\mapsto a_t(e)$ except for the placement of the ``initial'' jump time~$\tau_0$. For this we note that, as $t\mapsto a_t(e)$ is stationary, the random variable $U:=-\tau_0/(\tau_1-\tau_0)$ is uniform on~$[0,1]$ and independent of the family~\eqref{E:1.8a}. Starting from \eqref{E:1.8a} and an independent uniform~$U$, we just set $\tau_0:=-(T_0^\OFF+T_0^\ON)U$ and define the other~$\tau_i$ and~$\tau_i'$ accordingly. 
A minor complication is that the law of the interarrival times $\{(T_i^\OFF,T_i^\ON)\}_{i\in\Z}$ is not stationary under~$\BbbP$ but rather under the de-size-biased measure~$\wt\BbbP$ defined for $A\in\sigma(a_t(e)\colon t\in\R)$ by
\begin{equation}
\label{E:1.8}
\wt\BbbP(A):=\frac{\E((T_0^\OFF+T_0^\ON)^{-1}1_A)}{\E((T_0^\OFF+T_0^\ON)^{-1})},
\end{equation}
where, as before, $\E$ is expectation with respect to~$\BbbP$.

The representation based on \eqref{E:1.8a} and \eqref{E:1.8} makes it easier to describe specific examples. For instance, $\{(T_i^\OFF,T_i^\ON)\}_{i\in\Z}$ could be i.i.d.\ under~$\wt\BbbP$ which makes $t\mapsto a_t(e)$ a stationary renewal process modulo~2  under~$\BbbP$.  This is exactly the setting that many earlier studies (e.g., by Peres, Stauffer and Steif~\cite{PSS1}, Peres, Sousi and Steif~\cite{PSS2,PSS3} or Hermon and Sousi~\cite{HS}) have focused on. Another possibility is to draw $\{(T_i^\OFF,T_i^\ON)\}_{i\in\Z}$ from a stationary Markovian law on~$(0,\infty)\times(0,\infty)$ although even this is still too restrictive for our purposes.  Our result on dynamical percolation is thus cast as follows: 

\begin{theorem}
\label{dyn-perc}
Consider the random walk on dynamical percolation as specified above: The conductance processes $\{t\mapsto a_t(e)\}_{e\in E(\Z)}$ are i.i.d.\ taking values in $\{0,1\}$ with the associate sequence $\{(T_i^\OFF,T_i^\ON)\}_{i\in\Z}$ of interarrival times stationary  under~$\wt\BbbP$.  Assume, in addition to $T^\OFF_1,T^\ON_1$ being positive and finite, that
\begin{equation}
\label{E:1.9u}
\exists p>4\,\exists s>4\,\frac{1-1/p}{1-4/p}\colon\quad \wt\E\bigl((T_1^\OFF+T_1^\ON)^{p}\bigr)<\infty\,\wedge\,\wt\E\bigl((T_1^\ON)^{-s}\bigr)<\infty,
\end{equation}
where~$\wt\E$ is expectation with respect to~$\wt\BbbP$.
Then a Quenched Invariance Principle holds.
\end{theorem}

 The restriction to (at least) four moments of the ``OFF'' and ``ON'' times comes from that in Theorem~\ref{main}.  That being said,  some moment condition is  definitely needed  to ensure convergence to a non-degenerate Brownian motion. Indeed, as we show in Lemma~\ref{lemma-5.1}, when $T_i^\ON:=1$ for all~$i\in\Z$ and~$\{T_i^\OFF\}_{i\in\Z}$ are i.i.d.\ under~$\wt\BbbP$ with $\wt\E((T_1^\OFF)^{1/2})=\infty$, the random walk behaves subdiffusively. We do not know what moments of the ``ON/OFF''-times are critical for existence of such singular examples and/or the validity of a Quenched Invariance Principle.  We do not believe that our conditions \eqref{E:1.6} and \eqref{E:1.9u} are optimal; cf Remark~\ref{rem-non-optimal}. 

The specific example of random walk on dynamical percolation irregardless, the main thrust of our  result  is that it requires no assumptions (beyond stationarity and ergodicity  under space-time shifts)  on how the conductances evolve. This takes our approach significantly beyond earlier work (e.g., by B\'erard~\cite{B04}, Rassoul-Agha and Sepp\"alainen~\cite{RS05}, Bandyopadhyay and Zeitouni~\cite{BZ06}, Boldri\-ghini, Minlos and Pellegrinotti~\cite{BMP07}, Dolgopyat, Keller and Liverani~\cite{DKL08}, Redig and V\"ollering~\cite{RV13}) that requires more explicit assumptions. A limitation  of our approach compared to these studies is its  restriction to time-continuous variable speed random walks  with uniformly bounded jump rates. 

\section{Main steps and technical claims}
\label{sec-2}\noindent\nopagebreak
We proceed to discuss the main steps of the proof articulating the key technical statements to be established. The actual proofs come in Section~\ref{sec-3}.

\subsection{Overall picture}
There are two strategies we could follow in the proof of Theorem~\ref{main}. One could be based on elliptic regularity techniques developed earlier by Andres, Chiarini, Deuschel and Slowik~\cite{ACDS16} in $d\ge2$ and by Deuschel and Slowik~\cite{DS16} in $d=1$ for models satisfying, besides Assumption~\ref{ass-1}, suitable positive and negative moment conditions on the conductances. Besides the restriction to (a.s.) strictly positive conductances, a disadvantage of this approach is   the significant complexity caused by its reliance on  advanced techniques such as functional inequalities and Moser iteration.

The complexity notwithstanding, an important feature of the proofs in~\cite{ACDS16} and~\cite{DS16} is that the negative-moment condition is used only lightly --- mainly, to convert unadorned norms of important quantities to norms weighted by the conductances. In dimensions $d\ge2$, this was observed and fruitfully utilized by the first author and P.-F. Rodriguez~\cite{BR} to prove a Quenched Invariance Principle for models with bounded conductances assuming that the quantity in \eqref{E:1.4a}
obeys
\begin{equation}
\exists\varepsilon>0\colon\quad \E(T^{4d+\varepsilon})<\infty.
\end{equation}
While a similar (albeit still very technical) proof is expected to work for random walks with degenerate conductances in $d=1$, details of this have not been completed due to a different behavior of the Sobolev inequality in one spatial dimension.

Another  strategy we could follow would rely on  the aforementioned work of the first author~\cite{Bi1}. An inspection of the proofs of~\cite{Bi1} reveals that also here the negative moment condition is used only sporadically; namely, only in \cite[Lemma~4.2]{Bi1}, dealing with the construction of an auxiliary random walk that the whole proof is based on, and in \cite[Theorem~5.5]{Bi1} that constructs and proves the relevant properties of so called parabolic coordinates. We will follow this route and show that a slightly weaker form of Lemma~4.2 remains true, still sufficient to serve our purpose, and so does Theorem 5.5 provided we replace the negative-moment condition by assumptions (1-2) of Theorem~\ref{main}. 

\subsection{Main steps}
In order to bring the reader into the picture, let us recount the main steps of the proof in~\cite{Bi1}. The overall structure adheres to that of the proofs of invariance principles by the corrector method; see Biskup~\cite{B11} or Kumagai~\cite{Kumagai} for recent reviews. The proof thus starts with the construction of a \myemph{parabolic coordinate} which is a random map $\psi\colon\R\times\Z\to\R$ of which we require the following:
\settowidth{\leftmargini}{(1111)}
\begin{enumerate}
\item[(1)] $t\mapsto\psi(t,x)$ is continuous for each~$x\in\Z$ and $t,x\mapsto\psi(t,x)$ is a weak solution to
\begin{equation}
\label{E:2.1}
\frac\partial{\partial t}\psi(t,x)+L_t\psi(t,x)=0,\quad t\in\R,\,x\in\Z,
\end{equation}
with the ``initial'' data
\begin{equation}
\psi(0,0)=0.
\end{equation}
Here $L_t$ acts only on the second coordinate.
\item[(2)] For each~$t,s\in\R$ and each~$x,y\in\Z$, the cocycle condition holds
\begin{equation}
\label{E:2.5}
\psi(t+s,x+y)-\psi(t,x)=\psi(s,y)\circ\tau_{t,x}.
\end{equation}
\item[(3)] $\psi(\cdot,x)$ is, for each $x\in\Z$, a jointly measurable function of time (i.e., the first variable) and the random environment and we have
\begin{equation}
\label{E:2.6a}
\psi(t,x)\in L^1(\BbbP)\quad\text{\rm and}\quad\E\psi(t,x)=x,\quad t\in\R,\,x\in\Z,
\end{equation}
and
\begin{equation}
\label{E:2.7}
\E\bigl(a_0(0,1)\psi(0,1)^2\bigr)<\infty.
\end{equation}
\item[(4)] The spatial gradients of~$\psi(t,\cdot)$ are a.s.\ positive,
\begin{equation}
\label{E:2.8}
\psi(t,x+1)-\psi(t,x)>0,\quad t\in\R,\,x\in\Z.
\end{equation}
\end{enumerate}
Thinking of the map $x\mapsto\psi(t,x)$ as a different embedding of~$\Z$ into~$\R$, the above properties ensure that, in the new embedding, the random walk $t\mapsto\psi(t,X_t)$ is an $L^2$-martingale (under~$P^0_a$). 

Relying on the point of view of the particle enabled by Assumption~\ref{ass-1} and the Markov property of~$X$, we now check the conditions of the Functional Central Limit Theorem (see, e.g., Helland~\cite[Theorem~5.1(a)]{Helland}) for the process $t\mapsto\psi(t,X_t)$,  which thus tends in law, under a diffusive scaling of space and time, to Brownian motion with variance
\begin{equation}
\sigma^2:=2\E\bigl(a_0(0,1)\psi(0,1)^2\bigr).
\end{equation}
The proof of \cite[Theorem~1.2]{Bi1} contains all relevant (and explicit) details that apply to the present setting more or less \myemph{verbatim}.

While $\sigma^2<\infty$ by \eqref{E:2.7} and $\sigma^2>0$ is checked via \eqref{E:2.8}, the next, and usually the hardest, technical problem is to show that the ``deformation'' $\psi(t,X_t)-X_t$ of the random-walk path caused by the change of embedding of~$\Z$ is asymptotically irrelevant under the diffusive scaling of the process. As usual, it suffices to show this for the embedding itself which amounts to proving that the \myemph{parabolic corrector},
\begin{equation}
\label{E:2.8i}
\chi(t,x):=\psi(t,x)-x,
\end{equation}
obeys
\begin{equation}
\label{E:2.14}
\max_{\begin{subarray}{c}
x\in\Z\\|x|\le \sqrt n
\end{subarray}}
\,\,
\sup_{\begin{subarray}{c}
t\in\R\\0\le t\le n
\end{subarray}}
\frac{|\chi(t,x)|}{\sqrt n}
\,\underset{n\to\infty}\longrightarrow\,0,\quad\BbbP\text{\rm-a.s.}
\end{equation}
Indeed, the aforementioned Functional CLT gives $\max_{0\le t\le n}|\psi(t,X_t)|=O(\sqrt n)$ and \eqref{E:2.14} then shows $\max_{0\le t\le n}|\chi(t,X_t)|=o(\sqrt n)$  as desired. 

The proof of Theorem~\ref{main} is thus reduced to two technical steps: a construction of the parabolic coordinate~$\psi$ satisfying (1-4) above and a proof of the sublinear/subdiffusive bound \eqref{E:2.14}. In the approach of references \cite{ACDS16,DS16,BR}, this is exactly where elliptic regularity techniques are employed to their full extent. The approach of~\cite{Bi1} instead relies on the observation that, thanks to the one-dimensional nature of the problem, the spatial gradient of the parabolic coordinate
\begin{equation}
g(t,x):=\psi(t,x+1)-\psi(t,x)
\end{equation}
obeys the PDE
\begin{equation}
\label{E:2.11}
-\frac\partial{\partial t}g(t,x)=\mathcal{L}_t^+ g(t,x),
\end{equation}
where
the operator on the right-hand side acts on the spatial variable as
\begin{equation}
\mathcal{L}_{t}^{+}f(x):=b_{t}(x+1)f(x+1)+b_{t}(x-1)f(x-1)-2b_{t}(x)f(x)
\end{equation}
with
\begin{equation}
b_t(x):=a_t(x,x+1)
\end{equation} 
abbreviating the conductance of edge $(x,x+1)$.

As our use of adjoint notation suggests,~$\LL_t^+$ is the adjoint in $\ell^2(\Z)$ of an operator~$\LL_t$ that acts on $f\colon\Z\to\R$ as
\begin{equation}
\label{E:2.8a}
\LL_t f(x)= b_t(x)\bigl[\,f(x+1)+f(x-1)-2f(x)\bigr].
\end{equation}
A key point is that this is the generator of a continuous time simple symmetric random walk~$Y$ time-changed so that the jump rate at~$x$ at time~$t$ is~$2b_t(x)$,  which is a much simpler process than~$X$  to analyze. More importantly, the process~$Y$ also  provides all the needed tools for the proof of a Quenched Invariance Principle for the walk~$X$.

\subsection{Statements to be proved}
We will now describe what needs to be done in order to extend the proofs of \cite{Bi1} to that of Theorem~\ref{main}.
The first item of business is a formal construction of the random walk~$Y$. Note that this  walk moves on the set of edges of~$\Z$, which why we will refer to it as a \myemph{dual random walk}. The negative sign on the left of \eqref{E:2.11} necessitates that~$Y$ be run in negative time direction. The following generalizes
\cite[Lemma 4.2]{Bi1} to the situation when~$t\mapsto b_t(x)$ is allowed to vanish over sets of positive Lebesgue measure:

\begin{lemma}
\label{lemma-2.1}
Suppose that $t\rightarrow b_{-t}(x)$ is Borel-measurable and locally
integrable on $\left( 0,\infty \right) $ and,
in addition, for all $x\in \Z$,
\begin{equation}
\label{E:2.15}
\int_{0}^{\infty }b_{-t}(x)\,\textd t=\infty  
\end{equation}
Given~$x\in\Z$, let $P^x$ be the measure under which $Z$ is a discrete-time simple symmetric random walk on $\Z$ started from~$x$
and $N$ is an independent rate-1 Poisson point process. Then, for all $x\in \Z$, there is a non-decreasing continuous function $A\colon [0,\infty
)\rightarrow \lbrack 0,\infty )$ satisfying 
\begin{equation}
\label{E:2.16}
A(t)=\int_{0}^{t}2b_{-s}(Z_{N(A(s))})\textd s,\quad t\geq 0,
\end{equation}%
such that $P^{x}(A(t)<\infty )=1$ for each $t\geq 0$ and $x\in 
\Z$. Moreover, the process
\begin{equation}
\label{E:2.18ui}
Y_t:=Z_{N(A(t))},\quad t\ge0
\end{equation}
is a continuous-time Markov chain on~$\Z$ with generator~$\LL_t$ in \eqref{E:2.8a}.
\end{lemma}

Since $\BbbP$-a.s.\ validity of \eqref{E:2.15} is ensured by \eqref{E:1.6}, Lemma~\ref{lemma-2.1} shows that the dual random walk~$Y$ is well defined (as a time change of the  constant-speed  simple symmetric random walk) for $\BbbP$-a.e.\ sample of the random environment. 
As it turns out, the construction of the parabolic coordinate for~$X$ is equivalent to the construction of an invariant measure~$\Q$ for the environment as seen by the walk~$Y$. In~\cite{Bi1}, such an invariant measure is extracted by constructing directly its Radon-Nikodym derivative with respect to~$\BbbP$, and thus proving that~$\Q\ll\BbbP$. For us this comes in:

\begin{theorem}
\label{thm-2.2}
Under the conditions of Theorem~\ref{main}, there exists $\varphi \in L^{1}(\mathbb{P})$ that satisfies
\begin{enumerate}
\item[(1)] $\BbbP(\varphi>0)=1$ and $\mathbb{E}\varphi =1$,
\item[(2)] $\mathbb{E}(b_{0}(0)\varphi ^{2})\leq \mathbb{E}b_{0}(0)$,
\item[(3)] the map $t\rightarrow \varphi \circ \tau _{t,x}$ is
continuous and weakly differentiable such that 
\begin{equation}
\label{E:2.18}
\frac{\partial }{\partial t}\varphi \circ \tau _{t,x}+\mathcal{L}%
_{t}^{+}\varphi \circ \tau _{t,x}=0
\end{equation} 
holds for all~$t\in\R$ and~$x\in\Z$.
\end{enumerate}
In particular,~$\Q$ defined for~$A\in\FF$ by $\Q(A):=\E(\varphi 1_A)$ is a probability measure on $(\Omega,\FF)$ that is stationary and ergodic for the chain $t\mapsto\tau_{-t,Y_t}(a)$.
\end{theorem}

The link between the above Radon-Nikodym derivative and the parabolic coordinate is supplied by the observation that the PDEs \eqref{E:2.11} for $t,x\mapsto g(t,x)$ and \eqref{E:2.18} for $t,x\mapsto\varphi\circ\tau_{t,x}$ are identical. Setting $g(t,x):=\varphi\circ\tau_{t,x}$ would give us access to the gradient of~$\psi$. The parabolic coordinate~$\psi$ is extracted from this via
\begin{equation}
\label{E:5.2}
\psi(t,x):=\chi(t,0)\circ\tau_{0,x}+\sum_{k=0}^{x-1}\varphi\circ\tau_{0,k}
\end{equation}
where
\begin{equation}
\label{E:5.1}
\chi(t,0):=-\int_0^t\,\bigl(b_s(0)\varphi\circ\tau_{s,0}-b_s(-1)\varphi\circ\tau_{s,-1}\bigr)\,\textd s.
\end{equation}
Here the  (Lebesgue)  integral converges absolutely under expectation, and thus $\BbbP$-a.s.,\ by Tonelli's Theorem along with $\varphi\in L^1(\BbbP)$ and~$b_t(x)\in[0,1]$ as implied by the assumptions of Theorem~\ref{main}.
Standard interpretations of the integral in \eqref{E:5.1} and the sum in \eqref{E:5.2} are to be used for negative~$t$ and~$x$.\ 

It is straightforward to check (see the proof of \cite[Theorem~3.2]{Bi1}) that~$\psi$ from \eqref{E:5.2} obeys conditions (1-4) listed earlier. It then remains to prove the bound \eqref{E:2.14} for the corrector. (Note that \eqref{E:2.8i} and \eqref{E:5.1} are consistent.) In light of the cocycle conditions \eqref{E:2.5}, for this suffices to prove separately sublinearity in space
\begin{equation}
\label{E:6.4}
\lim_{n\to\pm\infty}\frac{|\chi(0,n)|}{n}=0,\quad\BbbP\text{\rm-a.s.}
\end{equation}
and subdiffusivity in time
\begin{equation}
\label{E:6.1}
\lim_{t\to\infty}\frac{|\chi(t,0)|}{\sqrt t}=0,\quad\BbbP\text{\rm-a.s.}
\end{equation}
Indeed, the ``good grid'' argument (originally designed in Berger and Biskup~\cite{BB07} for random walk on static percolation) used in~\cite{Bi1} then builds this into \eqref{E:2.14}.

As to the above almost sure limits, the one in \eqref{E:6.4} is proved by following the argument from \cite[Lemma~7.1]{Bi1} with the moment conditions supplied by Theorem~\ref{thm-2.2}(1). (In particular, no path interpolation as used in Berger and Biskup~\cite{BB07} are needed, nor is the conversion of the first moment of~$\chi$ to the weighted second moment from Biskup~\cite{B11}.) The proof of \eqref{E:6.1}, which comes as \cite[Proposition~7.2]{Bi1}, is considerably longer as it involves a different representation of $\chi(0,t)$ and the use of a Quenched Central Limit Theorem for the walk~$Y$ (which needs the invariant measure~$\Q$ and its equivalence with~$\BbbP$, as implied by Theorem~\ref{thm-2.2}(1)). But, as an inspection of these proofs reveals, the negative moment condition is not used throughout and some proofs (e.g., that of \cite[Lemma~8.1]{Bi1}) become even simpler for bounded conductances.

The bottom line is that the proof of Theorem~\ref{main} is reduced to those of Lemma~\ref{lemma-2.1} and Theorem~\ref{thm-2.2}. These proofs, which we will address in the next section, are the main technical contributions of the present note.



\section{Actual proofs}
\label{sec-3}
\noindent
We now move to the proofs of the technical claims from Section~\ref{sec-2}, starting with the construction of the dual random walk~$Y$ on which the rest of the argument is based. We assume the conditions of Theorem~\ref{main} throughout this section.

\subsection{The dual random walk}
As is standard in the theory of continuous-time Markov chains (cf., e.g., Liggett~\cite{Liggett-MC}), we first construct the transition probabilities of the desired Markov chain. For this we define a family of non-negative kernels $\cmss K_n(s,x;t,y)$ indexed by integers $n\ge0$ and depending on reals $-\infty<t\le s<\infty$ and vertices $x,y\in\Z$ inductively via 
\begin{multline}
\label{E:back-Kolm-n}
\qquad
\cmss K_{n+1}(s,x;t,y) := \texte^{-\int_t^s \,2b_u(x)\,\textd u}\,\delta_{x,y}
\\
+\int_t^s\,\texte^{-\int_r^s \,2b_u(x)\,\textd u}\,b_r(x)\,\Bigl(\sum_{z=\pm1}\cmss K_n(r,x+z;t,y)\Bigr)\,\textd r
\qquad
\end{multline} 
with the initial value $\cmss K_{0}(s,x;t,y):=0$. Note that time runs in the opposite direction of how the conductances are parametrized.

The definition \eqref{E:back-Kolm-n} readily yields that~$n\mapsto \cmss K_n(s,x;t,y)$ is non-decreasing and non-negative with $\sum_{y\in\Z}\cmss K_n(s,x;t,y)\le1$. The limit
\begin{equation}
\label{E:3.2}
\cmss K(s,x;t,y):=\lim_{n\to\infty}\cmss K_n(s,x;t,y)
\end{equation}
thus exists, is non-negative and obeys $\sum_{y\in\Z}\cmss K(s,x;t,y)\le1$ thanks to the Monotone Convergence Theorem. With all these objects being random variables on the probability space $(\Omega,\FF,\BbbP)$, we also have
\begin{equation}
\label{E:4.2ua}
\cmss K(s,x;t,y)\circ\tau_{u,z} = \cmss K(s+u,x+z;t+u,y+z)
\end{equation}
for all $s\ge t$, all~$u\in\R$ and all $x,y,z\in\Z$.  The next task is  to show that~$\cmss K$ is stochastic  which, as usual, is achieved by constructing the underlying Markov chain:  

\begin{proofsect}{Proof of Lemma~\ref{lemma-2.1}}
 As in the proof of \cite[Lemma 4.2]{Bi1}, instead of~$A$ we construct its inverse. Unfortunately, due to $s\mapsto b_{-s}(x)$ potentially vanishing over sets of positive Lebesgue measure, this inverse is no longer continuous which complicates its use. We thus proceed by a perturbation argument. 

Abusing our earlier notation, let~$\tau_0:=0<\tau_1<\dots$ denote the successive arrivals of a rate-1 (right-continuous) Poisson process~$N$ and let~$Z$ be the sample path of an independent discrete-time simple symmetric random walk on~$\Z$. Given~$\delta>0$, and restricting to the full-measure event~$\bigcap_{n\ge0}\{\tau_n<\infty\}$, set~$W_\delta(0):=0$ and, for each~$n\ge0$ and~$t\in(\tau_n,\tau_{n+1}]$, let~$W_\delta(t)$ be the unique number such that
\begin{equation}
\label{E:3.4i}
\int_{W_\delta(\tau_{n})}^{W_\delta(t)} \bigl[\delta+2b_{-s}(Z_{n})\bigr]\textd s=t-\tau
_{n}.
\end{equation}
The assumptions ensure that~$W_\delta(t)$ is finite for each~$t\ge0$ with~$t\mapsto W_\delta(t)$ continuous and strictly increasing with the lower bound $W_\delta(t)-W_\delta(s)\ge (2+\delta)^{-1}(t-s)$ whenever $t\ge s\ge0$. In particular, $\lim_{t\to\infty}W_\delta(t)=\infty$. 

It follows that $W_\delta$ admits a unique continuous and strictly increasing inverse~$A_\delta$ mapping $[0,\infty)$ onto itself. Thanks to the strict monotonicity, the  defining relation \eqref{E:3.4i} shows that $A_\delta(s)\in \lbrack\tau_n,\tau_{n+1})$ is equivalent to $s\in \lbrack W_\delta(\tau_n), W_\delta(\tau_{n+1}))$ and, since this forces $N(A_\delta(s))=n$ for all $s\in[W_\delta(\tau_n),W_\delta(\tau_{n+1}))$, we may rewrite \eqref{E:3.4i} into
\begin{equation}
\label{E:4.17ieu}
\int_{W_\delta(\tau_n)}^t \bigl[\delta+2b_{-s}(Z_{N(A_\delta(s))})\bigr]\textd s=A_\delta(t)-\tau_n,\quad t\in[W_\delta(\tau_n),W_\delta(\tau_{n+1})].
\end{equation} 
Here continuity of both sides in~$t$ was used to include $t=W_\delta(\tau_{n+1})$.

We now take~$\delta\downarrow0$ to extract the desired function~$A$. To that end we first note that a telescoping argument applied to \eqref{E:4.17ieu} gives
\begin{equation}
\label{E:3.6i}
\int_0^t\bigl[\delta+2b_{-s}(Z_{N(A_\delta(s))})\bigr]\textd s=A_\delta(t)
\end{equation}
for all~$t\ge0$, where we used that $A_\delta(0)=0$. A similar argument applied to \eqref{E:3.4i} shows that $\delta\mapsto W_\delta(t)$ is non-increasing and so $\delta\mapsto A_\delta(t)$ is non-decreasing. In light of the Lipschitz bound $0\le A_\delta(t)-A_\delta(s)\le (2+\delta)(t-s)$ for $t\ge s\ge0$, the limit
\begin{equation}
A(t):=\lim_{\delta\downarrow0}A_\delta(t)
\end{equation}
 exists and defines a continuous  real-valued  non-decreasing function $t\mapsto A(t)$ with~$A(0)=0$. The upward monotonicity of $\delta\mapsto A_\delta(s)$ in conjunction with the  right-continuity of~$N$ gives $b_{-s}(Z_{N(A_\delta(s))})\to b_{-s}(Z_{N(A(s))})$ as~$\delta\downarrow0$, for each~$s\ge0$. Taking $\delta\downarrow0$ in \eqref{E:3.6i} with the help of the Bounded Convergence Theorem then proves \eqref{E:2.16}. 

Define~$Y$ from the processes~$N$, $Z$ and~$A$ by the formula \eqref{E:2.18ui}. Recall that~$P^x$ is the law of these objects such that $P^x(Z_0=x)=1$. The identity \eqref{E:back-Kolm-n} then inductively  shows that, for all~$t\ge0$,
\begin{equation}
\label{E:3.16}
\cmss K_n(0,x;-t,y)=P^x\bigl(Y_t=y,\,N(A(t))<n\bigr),\quad n\ge0.
\end{equation}
As~$P^x(N(A(t))<\infty)=1$ due to $P^x(A(t)<\infty)=1$, taking~$n\to\infty$ gives
\begin{equation}
\cmss K(0,x;-t,y)=P^x(Y_t=y).
\end{equation}
In particular,~$\cmss K$ is stochastic and, taking~$n\to\infty$ in \eqref{E:back-Kolm-n} using the Monotone Convergence Theorem,~$Y$ is a Markov chain with generator~$\LL_t$.
\end{proofsect}

\subsection{Proof of Theorem~\ref{thm-2.2}}
Having constructed the random walk~$Y$, we now move to the construction of the Radon-Nikodym derivative~$\varphi$ of the invariant measure on environments as seen from~$Y$. 
 As in~\cite{Bi1}, we  will extract~$\varphi$ as an $\varepsilon\downarrow0$ limit of the quantity 
\begin{equation}
\label{E:3.9}
\varphi _{\varepsilon }:=\varepsilon \int_{0}^{\infty }\texte^ {-\varepsilon
t}\Bigl(\,\sum_{y\in \Z}\cmss K(t,y;0,0)\Bigr)\textd t.
\end{equation}%
We first pull some  observations from \cite{Bi1}:

\begin{lemma}
\label{heat}
For each~$\varepsilon>0$,
\begin{equation}
\label{E:3.10}
\mathbb{E}\varphi _{\varepsilon }=1
\end{equation}
and, in particular, $\varphi_\varepsilon\in[0,\infty)$ $\BbbP$-a.s. Moreover, abbreviating
\begin{equation}
\varphi _{\varepsilon }(t,x):=\varphi _{\varepsilon }\circ \tau
_{t,x}
\end{equation}
the function $t\mapsto \varphi _{\varepsilon }(t,x)$ is continuous and weakly differentiable with
\begin{equation}
\label{E:3.12}
\frac{\partial }{\partial t}\varphi _{\varepsilon }(t,x)=\varepsilon \left(
\varphi _{\varepsilon }(t,x)-1\right) -\mathcal{L}_{t}^{+}\varphi
_{\varepsilon }(t,x).
\end{equation}
\end{lemma}

\begin{proofsect}{Proof}
Formula \eqref{E:3.10} is obtained by invoking stationarity of~$\BbbP$ along with \eqref{E:4.2ua} and the Monotone Convergence Theorem to rewrite the sum in \eqref{E:3.9} under expectation into $\sum_{y\in\Z}\cmss K(0,0;-t,y)=1$. Formula \eqref{E:3.12} is a limit version of \cite[formula~5.20]{Bi1} whose derivation applies \myemph{verbatim}.
\end{proofsect}

\begin{lemma}
\label{lemma-3.2}
For each $\varepsilon >0$, we
have%
\begin{equation}
\label{E:3.13}
\mathbb{E}\bigl(b_{0}(0)\varphi _{\varepsilon }^{2}\bigr)\leq \mathbb{E}b_{0}(0).
\end{equation}
\end{lemma}

\begin{proofsect}{Proof}
This is a restatement of \cite[Proposition 5.1]{Bi1} whose proof applies without changes in our case as well.
\end{proofsect}

The argument of~\cite{Bi1} proceeds by taking a weak limit of~$\varphi_\varepsilon$ as~$\varepsilon\downarrow0$ and using  \eqref{E:3.13}  to show that ``no mass is lost'' in \eqref{E:3.10} in this process. In~\cite{Bi1}, this step required the negative moment condition which restricts us to $b_0(0)>0$ $\BbbP$-a.s. Once this does not apply, even the  subsequent  use of  \eqref{E:3.13}  becomes problematic as the inequality can at best provide control of the weak limit only on the set where $b_0(0)>0$.  

In order to overcome these issues, we invoke an idea from Biskup and Rodriguez~\cite{BR} that is itself drawn from Mourrat and Otto~\cite{MO}. In these works, the argument proceeds by finding a version of  \eqref{E:3.13}  in which~$b_t(x)$ is replaced by the time-averaged quantity of the form
\begin{equation}
c_{t}(x):=\int_{t}^{\infty }k_{s-t}b_{s}(x)\,\textd s,
\end{equation}%
where~$t\mapsto k_t$ is a suitable positive function on~$(0,\infty)$ with sufficient decay at infinity. Note that $c_t(x)$ is positive as soon as~$s\mapsto b_s(x)$ is positive on a set of positive Lebesgue measure, which for us occurs $\BbbP$-a.s.\ thanks to \eqref{E:1.6}. 

As it turns out, the most useful choice is to take~$t\mapsto k_t$ with a power-law decay and so we henceforth set
\begin{equation}
\label{E:3.15}
k_{t}:=(1+t)^{-\alpha }
\end{equation}
for some $\alpha>0$ to be determined momentarily. The reason for this is seen from:

\begin{lemma}
\label{lemma-3.3}
Recall the quantity~$T$ from \eqref{E:1.4a} and let $\alpha>0$. Then for~$c_0(0)$ defined using the kernel \eqref{E:3.15},
\begin{equation}
\E(T^{\alpha })<\infty\quad\Rightarrow\quad\E(c_0(0)^{-1})<\infty.
\end{equation}
\end{lemma}

\begin{proofsect}{Proof}
We have 
\begin{equation}
c_{0}(0) =\int_{0}^{\infty }k_{t}b_{t}(0)\,\textd t \ge \int_{0}^{T}k_{t}b_{t}(0)\,\textd t 
\ge\frac{1}{(1+T)^{\alpha }}\int_{0}^{T}b_{t}(0)\textd t=\frac{1}{(1+T)^{\alpha }},
\end{equation}
which yields $c_{0}(0)^{-1}\le (1+T)^{\alpha}$. The claim follows.
\end{proofsect}

The restriction on the moment of~$T$ now comes via:

\begin{proposition}
\label{prop-3.4}
For any~$\alpha>3$ and for~$c_0(0)$ defined using the kernel \eqref{E:3.15},
\begin{equation}
\label{E:3.18}
\sup_{0<\varepsilon<1}\E\bigl(c_0(0)\varphi_\varepsilon^2\bigr)<\infty.
\end{equation}
\end{proposition}

Before giving the proof of Proposition~\ref{prop-3.4}, which comes in Section~\ref{sec-3.3}, we give:

\begin{proofsect}{Proof of Theorem~\ref{thm-2.2} from Proposition~\ref{prop-3.4}}
Suppose the moment condition \eqref{E:1.6} holds with some~$\varepsilon>0$ and let~$\alpha:=3+\varepsilon$. Writing $L^0(\BbbP)$ for the set of measurable~$f\colon\Omega\to\R$ modulo changes on $\BbbP$-null sets, consider the Hilbert space
\begin{equation}
\HH:=\Bigl\{f\in L^0(\BbbP)\colon \E\bigl(c_0(0)^{-1}f^2\bigr)<\infty\Bigr\}
\end{equation}
endowed with the inner product $\langle f,g\rangle_\HH:=\E(c_0(0)^{-1}fg)$.
 Using~$C$ to denote  the supremum in \eqref{E:3.18}, for any $f\in L^{\infty }(\mathbb{P})$ and $\varepsilon\in(0,1)$ the Cauchy-Schwarz inequality shows
\begin{equation}
\mathbb{E}(\varphi _{\varepsilon }f)\le C^{1/2}\bigl[ \mathbb{E}(c_{0}(0)^{-1}f^{2})\bigr] ^{1/2}.
\end{equation}
It follows that
\begin{equation}
\phi_\varepsilon(f):= \mathbb{E}(\varphi _{\varepsilon }f)
\end{equation}
 defines a continuous
linear functional on~$\HH$ with the operator norm bounded
by $C^{1/2}$ uniformly in~$\varepsilon\in(0,1)$. As  $\HH$ is separable, and the unit ball in~$\HH^\star$ thus weakly compact, the Cantor diagonal argument yields a sequence $\varepsilon_n\downarrow0$ and~$\phi\in\HH^\star$ such that $\phi_{\varepsilon_n}(f)\to\phi(f)$ for all~$f\in\HH$.  The Riesz lemma then shows that~$\phi$ takes the form $\phi(f)=\mathbb{E}%
\left[ c_{0}(0)^{-1}hf\right] $ for some $h\in \HH$.  We define  $\varphi :=c_{0}(0)^{-1}h$.

Lemma~\ref{lemma-3.3}  along with the moment condition \eqref{E:1.6} implies  that the space $L^\infty(\BbbP)$ of bounded measurable functions obeys 
\begin{equation}
L^{\infty }(\mathbb{P})\subset \HH.
\end{equation}
In particular, $1\in\HH$. The identity $\phi_\epsilon(1)=\mathbb{E}\varphi _{\varepsilon }=1$ then survives the limit and so we get $\E\varphi=1$, proving the second half of~(1). For the inequality in~(2), we first note that the bound $\mathbb{E}[b_{0}(0)(\varphi _{\varepsilon
}-f)^{2}]\ge0$ shows that, for any $f\in L^\infty(\BbbP)$,
\begin{equation}
2\mathbb{E}(b_{0}(0)f\varphi
_{\varepsilon })-\mathbb{E}(b_{0}(0)f^{2})\le\mathbb{E}[b_{0}(0)\varphi
_{\varepsilon }^{2}]\leq \mathbb{E}[b_{0}(0)],
\end{equation}%
where the last inequality is taken from Lemma~\ref{lemma-3.2}. Since $b_{0}(0)f\in
L^{\infty }(\mathbb{P})$ for~$f\in L^\infty(\BbbP)$, the first term on the left converges to $2\E(b_0(0)f\varphi)$ along the sequence $\{\varepsilon_n\}_{n\ge1}$  that was used to define~$\varphi$.
Combining this with
\begin{equation}
\mathbb{E}[b_{0}(0)\varphi ^{2}]=\sup_{f\in L^{\infty }(\mathbb{P})}\bigl[2%
\mathbb{E}(b_{0}(0)f\varphi )-\mathbb{E}(b_{0}(0)f^{2})\bigr]
\end{equation}
which, as is checked by a suitable truncation, holds regardless whether the left-hand side is finite or infinite, then yields the inequality in~(2).

The proof of the needed regularity of $t\mapsto\varphi\circ\tau_{t,x}$ --- or, more precisely, the existence of a continuous, weakly-differentiable version --- so that the PDE \eqref{E:2.18} holds is identical to that in \cite{Bi1} and we omit it here. It remains to prove the $\BbbP$-a.s.\ positivity of~$\varphi$. First note that $\varphi$ is
nonnegative. This is because $%
\mathbb{E}(1_{\{\varphi<0\}}\varphi _{\varepsilon })$ tends to $\mathbb{E}(1_{\{\varphi<0\}}\varphi)$ in the limit defining~$\varphi$ and $\varphi _{\varepsilon }\geq 0$ then shows $\mathbb{E}(1_{\{\varphi<0\}}\varphi)\ge0$ forcing~$\BbbP(\varphi<0)=0$. Next we observe that, for each~$t\ge0$ we have
\begin{equation}
\label{E:3.24}
\varphi = \sum_{x\in\Z}\varphi\circ\tau_{t,x}\,\cmss K(t,x;0,0)
\end{equation}
on a set of full~$\BbbP$-measure, which is proved using the same argument as in~\cite{Bi1}. As $\cmss K(t,0;0,0)>0$, assuming~$\varphi=0$ in \eqref{E:3.24} forces $\varphi\circ\tau_{t,0}=0$ $\BbbP$-a.s.\ for each~$t\ge0$. Using shift invariance and continuity, we conclude
\begin{equation}
\label{E:3.25}
\{\varphi=0\}\,\overset{\BbbP\text{-a.s.}}=\,\bigl\{\forall t\in\R\colon \varphi\circ\tau_{t,0}=0\bigr\}.
\end{equation}
But for each~$x\in\Z$ and $\BbbP$-a.e.\ realization of the random environment, $\cmss K(t,x;0,0)>0$ once~$t\ge0$ is large enough and so, invoking \eqref{E:3.24} and shift invariance again we get
\begin{equation}
\label{E:3.26}
\{\varphi=0\}\,\overset{\BbbP\text{-a.s.}}=\,\bigl\{\forall t\in\R\,\forall x\in\Z\colon \varphi\circ\tau_{t,x}=0\bigr\}.
\end{equation}
The event on the right is shift invariant and so, in light of ergodicity of~$\BbbP$ from Assumption~\ref{ass-1}, it is a zero-one event under~$\BbbP$. The case of full measure is ruled out by~$\E\varphi=1$ thus proving $\BbbP(\varphi=0)=0$.
\end{proofsect}

\subsection{Boundedness of weighted Dirichlet energy}
\label{sec-3.3}\noindent
The last remaining item needed to complete the proof of Theorem~\ref{main} is the proof of the uniform bound \eqref{E:3.18}. We again need a couple lemmas that are drawn from, or otherwise available in~\cite{Bi1}.
Define 
\begin{equation}
\chi _{\varepsilon }:=\int_{0}^{\infty }\texte^{-\varepsilon t}\bigl[
b_{t}(0)\varphi _{\varepsilon }(t,0)-b_{t}(-1)\varphi _{\varepsilon }(t,-1)%
\bigr] \textd t,
\end{equation}%
where the integral of each of the two terms in the square bracket is finite under expectation with respect to~$\BbbP$, and thus~$\BbbP$-a.s., by the fact that $0\le b_t(0)\varphi_\varepsilon(t,0)\le \varphi_\varepsilon(t,0)$ thanks to \eqref{E:1.4a} and $\varphi_\varepsilon(t,0)\in L^{1}(\mathbb{P})$ thanks to \eqref{E:3.10}. We start with:

\begin{lemma}
\label{L2chi}For each $\varepsilon >0$, 
\begin{equation}
\left\Vert \chi _{\varepsilon }\right\Vert _{L^{2}(\mathbb{P})}\leq \frac{2}{%
\varepsilon }.
\end{equation}
\end{lemma}

\begin{proofsect}{Proof}
Minkowski's inequality yields
\begin{equation}
\begin{aligned}
\Vert \chi _{\varepsilon }\Vert _{L^{2}(\mathbb{P})} 
&\leq \int \texte^{-\varepsilon t}\bigl\Vert b_{t}(0)\varphi _{\varepsilon
}(t,0)-b_{t}(-1)\varphi _{\varepsilon }(t,-1)\bigr\Vert _{L^{2}(\mathbb{P}%
)}\textd t \\
&\leq 2\int \texte^{-\varepsilon t}\bigl\Vert b_{t}(0)\varphi _{\varepsilon
}(t,0)\bigr\Vert _{L^{2}(\mathbb{P})}\textd t 
\leq 2\int \texte^{-\varepsilon t}\bigl\Vert b_{t}(0)^{1/2}\varphi
_{\varepsilon }(t,0)\bigr\Vert _{L^{2}(\mathbb{P})}\textd t,
\end{aligned}
\end{equation}
where we used $b_{t}(0)\leq 1$ in the last inequality. Lemma~\ref{lemma-3.2} along with~$\E b_0(0)\le1$ bounds the last $L^2$-norm by one.
\end{proofsect}

The motivation for introducing $\chi_\epsilon$ in \cite{Bi1} is that its spatial gradients (under environment shifts) are those of centered~$\varphi_\epsilon$, which  (in light of~$\varphi$ being the gradient of the parabolic coordinate)  makes~$\chi_\epsilon$ an approximate corrector. Indeed, we have:

\begin{lemma}
\label{lemma 5.2}
For each $\varepsilon >0$, 
\begin{equation}
\label{E:3.31}
\chi _{\varepsilon }\circ \tau _{0,1}-\chi _{\varepsilon }=\varphi
_{\varepsilon }-1.
\end{equation}
\end{lemma}

\begin{proofsect}{Proof}
\cite[Lemma~5.2]{Bi1} proves a truncated version of this equation; namely,%
\begin{equation}
\label{E:3.32}
\chi _{\varepsilon ,n}\circ \tau _{0,1}-\chi _{\varepsilon ,n}=\varphi
_{\varepsilon ,n+1}-1,
\end{equation}%
where $\varphi _{\varepsilon ,n}$ is defined by \eqref{E:3.9} with~$\cmss K$ replaced by~$\cmss K_n$ and 
\begin{equation}
\chi _{\varepsilon ,n}:=\int_{0}^{\infty }\texte^{-\varepsilon t}\bigl[
b_{t}(0)\varphi _{\varepsilon ,n}\circ \tau _{t,0}-b_{t}(-1)\varphi
_{\varepsilon ,n}\circ \tau _{t,-1}\bigr] \textd t.
\end{equation}%
The monotonicity of~$n\mapsto\cmss K_n$ implies $\varphi _{\varepsilon ,n}\uparrow \varphi _{\varepsilon }$ as $%
n\rightarrow \infty $ and the Monotone Convergence Theorem shows $\chi _{\varepsilon ,n}\to
\chi _{\varepsilon }$ as~$n\to\infty$. Hence \eqref{E:3.31} follows from \eqref{E:3.32}.
\end{proofsect}

Lemma~\ref{lemma 5.2} now extends the bound from Lemma~\ref{L2chi} to~$\varphi_\varepsilon$ as well:

\begin{lemma}
\label{L2phi} 
For each~$\varepsilon>0$,
\begin{equation}
\mathbb{E}\varphi _{\varepsilon }^{2}\leq \Bigl(1+\frac4\varepsilon\Bigr)^2.
\end{equation}
\end{lemma}

\begin{proofsect}{Proof}
The triangle inequality applied to the identity from Lemma~\ref{lemma 5.2} gives
\begin{equation}
\Vert \varphi _{\varepsilon }\Vert _{L^{2}(\mathbb{P})} 
\leq
1+\Vert \chi _{\varepsilon }\circ \tau _{0,1}\Vert _{L^{2}(%
\mathbb{P})}+\Vert \chi _{\varepsilon }\Vert _{L^{2}(\mathbb{P})}
= 1+2\Vert \chi _{\varepsilon }\Vert _{L^{2}(\mathbb{P})}.
\end{equation}
Lemma~\ref{L2chi} now bounds the right-hand side by $1+\frac4\varepsilon$.
\end{proofsect}

With the above lemmas in hand, we are ready to give:

\begin{proofsect}{Proof of Proposition~\ref{prop-3.4}}
Our task is to convert the Dirichlet energy with averaged conductance to the Dirichlet energy with instantaneous conductance to which the inequality in Lemma~\ref{lemma-3.2} can be applied. As observed first in Mourrat and Otto~\cite{MO} and further exploited in Biskup and Rodriguez~\cite{BR}, this  is possible thanks to the fact that  $t,x\mapsto\varphi_\varepsilon(t,x)$  obeys  the (massive) heat equation \eqref{E:3.12}. We start with the rewrite
\begin{equation}
\label{E:3.36}
\begin{aligned}
c_{0}(0)\varphi _{\varepsilon }(0,0)^{2} 
&=\int_{0}^{\infty }k_{t}b_{t}(0)\varphi _{\varepsilon }(0,0)^{2}\textd t \\
&=\int_{0}^{\infty }k_{t}b_{t}(0)\bigl[ \varphi _{\varepsilon
}(0,0)-\varphi _{\varepsilon }(t,0)+\varphi _{\varepsilon }(t,0)\bigr]
^{2}\textd t  
\\
&\leq 2\int_{0}^{\infty }k_{t}b_{t}(0)\varphi _{\varepsilon
}(t,0)^{2}\textd t+2\int_{0}^{\infty }k_{t}b_{t}(0)\bigl( \varphi _{\varepsilon
}(0,0)-\varphi _{\varepsilon }(t,0)\bigr) ^{2}\textd t  
\\
&\leq 2\int_{0}^{\infty }k_{t}b_{t}(0)\varphi _{\varepsilon
}(t,0)^{2}\textd t+2\int_{0}^{\infty }k_{t}\bigl( \varphi _{\varepsilon
}(0,0)-\varphi _{\varepsilon }(t,0)\bigr) ^{2}\textd t,
\end{aligned}
\end{equation}
where we use the inequality $(a+b)^{2}\leq 2a^{2}+2b^{2}$ and the assumption 
$b_{t}(0)\leq 1$. For the integrand of the second term, the heat
equation in Lemma~\ref{heat} along with the Cauchy-Schwarz inequality and the bound $(\sum_{i=1}^4 a_i)^2\le 4\sum_{i=1}^4 a_i^2$ yields
\begin{equation}
\begin{aligned}
\bigl( \varphi _{\varepsilon }&(0,0)-\varphi _{\varepsilon }(t,0)\bigr)
^{2} 
=\left[ \int_{0}^{t}\mathcal{L}_{s}^{+}\varphi _{\varepsilon
}(s,0)-\varepsilon \left( \varphi _{\varepsilon }(t,x)-1\right) \textd s\right]
^{2} \\
&=\Bigl[ \int_{0}^{t}b_{s}(1)\varphi _{\varepsilon }(s,1)+b_{s}(-1)\varphi
_{\varepsilon }(s,-1)-2b_{s}(0)\varphi _{\varepsilon }(s,0)-\varepsilon
\left( \varphi _{\varepsilon }(s,0)-1\right) \textd s\Bigr] ^{2} \\
&\leq t\int_{0}^{t}\bigl[ b_{s}(1)\varphi _{\varepsilon
}(s,1)+b_{s}(-1)\varphi _{\varepsilon }(s,-1)-2b_{s}(0)\varphi _{\varepsilon
}(s,0)-\varepsilon \left( \varphi _{\varepsilon }(s,0)-1\right) \bigr]
^{2}\textd s \\
&\leq 4t\varepsilon ^{2}\int_{0}^{t}( \varphi _{\varepsilon
}(s,0)-1) ^{2}\textd s+16t\sum_{z=-1}^1\int_{0}^{t}\bigl[ b_{s}(z)\varphi _{\varepsilon
}(s,z)\bigr]^2\textd s,
\end{aligned}
\end{equation}
where the factor 16 involves overcounting that makes the resulting expression simpler to write. Bounding $(\varphi_\varepsilon(s,0)-1)^2\le 2+2\varphi_\varepsilon(s,0)^2$ and   using $b_s(z)\le1$ to drop one $b_s(z)$ from  the second integral wraps this into
\begin{equation}
\bigl( \varphi _{\varepsilon }(0,0)-\varphi _{\varepsilon }(t,0)\bigr)
^{2} 
\le
8t^2\varepsilon^2+8t\varepsilon^2\int_0^t \varphi_\varepsilon(s,0)^2\textd s+
16t\sum_{z=-1}^1\int_{0}^{t}b_{s}(z)\varphi _{\varepsilon
}(s,z)^2\textd s.
\end{equation}
Plugging the resulting bound on the right of \eqref{E:3.36} and performing a simple change of variables then shows
\begin{multline}
\quad
c_{0}(0)\varphi _{\varepsilon }(0,0)^{2} 
\le
2\int_0^\infty k_tb_t(0)\phi_\epsilon(t,0)^2\textd t+16\varepsilon^2\int_0^\infty t^2k_t\,\textd t
\\+16\varepsilon^2\int_0^\infty K_t \varphi_\varepsilon(t,0)^2\textd t
+
32t\sum_{z=-1}^1\int_{0}^{\infty}K_t b_{t}(z)\varphi _{\varepsilon
}(t,z)^2\textd t,
\quad
\end{multline}
where
\begin{equation}
K_t:=t\int_t^\infty k_s\,\textd s.
\end{equation}
Taking expectation and invoking stationarity of~$\BbbP$ with respect to shifts gives
\begin{multline}
\qquad
\E\bigl(c_{0}(0)\varphi _{\varepsilon }(0,0)^{2} \bigr)
\le 2\E\bigl(b_0(0)\varphi_\varepsilon^2\bigr)\Bigl(\,\int_0^\infty k_t\,\textd t\Bigr)
 + 16\varepsilon^2\Bigl(\,\int_0^\infty t^2k_t\,\textd t \Bigr)
\\+ \Bigl(\int_0^\infty K_t\,\textd t\Bigr)\Bigl[16\varepsilon^2\E(\varphi_\varepsilon^2)+96\E\bigl(b_0(0)\varphi_\varepsilon^2\bigr)\Bigr].
\qquad
\end{multline}
For our choice \eqref{E:3.15} with~$\alpha>3$, the integrals with respect to~$t$ converge and the terms involving the expectations are bounded uniformly in~$\varepsilon\in(0,1)$ thanks to Lemma~\ref{L2phi} and Theorem~\ref{thm-2.2}(2).
\end{proofsect}

\begin{remark}
\label{rem-non-optimal}
Similarly as in the derivations of~\cite{BR}, the use of Cauchy-Schwarz inequality along with dropping factors of $b_s(\cdot)$ is likely a wasteful step that forces the need for higher moments of~$T$ than what should be optimal and limits us to bounded conductances. However, we do not know how to proceed otherwise. 
\end{remark}

\section{Random walk on dynamical percolation}
\noindent
We will now apply the conclusions of Theorem~\ref{main} to random walk on dynamical percolation. Recall that, in our interpretation, a dynamical percolation is any conductance environment  with law~$\BbbP$ under  which  $\{t\mapsto a_t(e)\}_{e\in E(\Z)}$ are  i.i.d.\ copies of a  zero-one valued, non-degenerate (i.e., truly two-valued), stationary  process on~$\{0,1\}$ with piece-wise constant right-continuous sample paths. A standard argument gives:

\begin{lemma}
Any dynamical percolation law~$\BbbP$ obeys Assumption~\ref{ass-1}.
\end{lemma}

\begin{proofsect}{Proof}
The required regularity of sample paths follows from the assumed piece-wise constancy. The law~$\BbbP$ is also clearly invariant under all space time shifts. In order to show ergodicity, let~$A\in\FF$ be invariant under the space shifts (invariance under time shifts is not required). The product structure of~$A$ ensures that, given $n\ge1$, there exists $A_n\in\sigma(a_t(x,x+1)\colon -n\le x<n)$ such that $\E|1_A-1_{A_n}|<1/n$. Now define $A_n':=\tau_{0,n}(A_n)$ and $A_n'':=\tau_{0,-n}(A_n)$ and use space-shift invariance of~$A$ to check the inequalities $\BbbP(A_n')<\BbbP(A)+1/n$, $\BbbP(A_n'')<\BbbP(A)+1/n$ and $\BbbP(A)-\BbbP(A_n'\cap A_n'')<2/n$. Observing that $A_n'$ and~$A_n''$ are independent under~$\BbbP$ and taking~$n\to\infty$ this yields $\BbbP(A)\le\BbbP(A)^2$, thus showing that~$A$ is trivial under~$\BbbP$.
\end{proofsect}

With this in hand, we now give:

\begin{proofsect}{Proof of Theorem~\ref{dyn-perc}}
Thanks to the assumed non-degeneracy, we can identify each individual conductance process $t\mapsto a_t(e)$ with a sequence $\{(T_i^\OFF,T_i^\OFF)\}_{i\in\Z}$ of positive and finite random variables which is stationary under the law~$\wt\BbbP$. The restriction of the law~$\BbbP$ to $A\in\sigma(\{a_t(e)\colon t\in \R\})$ is then obtained by inverting \eqref{E:1.8} to
\begin{equation}
\label{E:4.1a}
\BbbP(A)=\frac{\wt\E\bigl((T^\OFF_0+T^\ON_0)1_A\bigr)}{\wt\E(T^\OFF_0+T^\ON_0)},
\end{equation}
where the expectations exist thanks to our moment assumption in \eqref{E:1.9u}. (The law~$\BbbP$ on  full environment  requires one size-biasing factor for each edge.)

Using the ``OFF/ON''-times, the random variable~$T$ from \eqref{E:1.4a} can be bounded as
\begin{equation}
\label{E:4.1}
T\le \sum_{i=0}^N(T_i^\OFF+T_i^\ON),
\end{equation}
where
\begin{equation}
\label{E:4.2}
N:=\inf\Bigl\{n\ge1\colon\sum_{i=1}^n T_i^\ON\ge1\Bigr\}.
\end{equation}
In order to estimate the moments of the sum on the right of \eqref{E:4.1}, we recall the following observation from Berger and Biskup~\cite{BB07}:

\begin{lemma}[Lemma~4.5 of~\cite{BB07}]
\label{lemma-4.1}
Given reals $p>1$, $r\in[1,p)$ and~$s$ such that
\begin{equation}
\label{E:4.3}
s>r\frac{1-1/p}{1-r/p},
\end{equation}
if~$X_1,X_2,\dots$ are random variables such that $\sup_{i\ge1}\Vert X_i\Vert_p<\infty$ and~$N$ is integer valued such that~$N\in L^s$, then $\sum_{i=1}^N X_i\in L^r$.
\end{lemma}

In order to apply this to our situation, let~$p$ and $s$ be reals satisfying  the inequalities and  the moment bounds in \eqref{E:1.9u}.  Pick any~$\tilde s$ satisfying \begin{equation}
s>\tilde s>4\frac{1-1/p}{1-4/p}.
\end{equation}
  Continuity then ensures that there is $r\in(4,p)$ such that \eqref{E:4.3} holds  with $\tilde s$ in place of~$s$. 
In order to apply Lemma~\ref{lemma-4.1}, we need to control the moments of~$N$ in \eqref{E:4.2}. Here the Markov and Jensen inequalities show
\begin{equation}
\begin{aligned}
\wt\BbbP(N>n)&\le \wt\BbbP\biggl(\,\sum_{i=1}^n T_i^\ON<1\biggr)
\le \wt\E\biggl(\Bigl(\,\sum_{i=1}^n T_i^\ON\Bigr)^{-s}\biggr)
\\&=\frac1{n^s} \wt\E\biggl(\Bigl(\frac1n\sum_{i=1}^n T_i^\ON\Bigr)^{-s}\biggr)
\le\frac1{n^s}\wt\E\biggl(\frac1n\sum_{i=1}^n (T_i^\ON)^{-s}\biggr)=\frac1{n^s}\wt\E\bigl((T_1^\ON)^{-s}\bigr)
\end{aligned}
\end{equation}
and  the formula $\wt\E(N^{\tilde s})=\int_0^\infty \tilde s n^{\tilde s-1}\wt\BbbP(N>n)\textd n$ then gives  $N\in L^{\tilde s}$.  Lemma~\ref{lemma-4.1}  (with $\tilde s$ in place of~$s$)  then shows $\wt\E(T^r)<\infty$. In order to convert this to a bound under expectation with respect to~$\BbbP$, we invoke the H\"older inequality to get
\begin{equation}
\E(T^\alpha)= 
\frac{\wt\E\bigl((T^\OFF_0+T^\ON_0)T^\alpha\bigr)}{\wt\E(T^\OFF_0+T^\ON_0)}
\le\frac{\bigl[\wt\E((T^\OFF_0+T^\ON_0)^4)\bigr]^{1/4}}{\wt\E(T^\OFF_0+T^\ON_0)}\bigl[\wt\E(T^{4\alpha/3})\bigr]^{3/4}.
\end{equation}
Note that the fourth moment exists by our assumption in \eqref{E:1.9u}.
Setting $\alpha:=3r/4$ and noting that then $\alpha>3$, we have verified the moment condition \eqref{E:1.6}. Theorem~\ref{main} shows that a Quenched Invariance Principle holds.
\end{proofsect}

\begin{lemma}
\label{lemma-5.1}
Suppose that the conductances are independent with the associated ``OFF/ON''-times such that $T_i^\ON:=1$ for all~$i\in\Z$ and~$\{T_i^\OFF\}_{i\in\Z}$ i.i.d.\ under~$\wt\BbbP$ with   $T_0^\OFF\not\in L^{1/2}(\wt\BbbP)$.  Then $X_t/\sqrt t\to0$ in probability as $t\to\infty$ for~$\BbbP$-a.e.\ sample of the random environment.
\end{lemma}

\begin{proofsect}{Proof}
For each edge~$e\in E(\Z)$, let
\begin{equation}
\wt T(e):=\inf\bigl\{t\ge0\colon a_t(e)>0\bigr\}.
\end{equation}
Under the assumptions of the lemma, and with the size-biasing in \eqref{E:4.1a} taken into account, $\{\wt T(e)\}_{e\in E(\Z)}$ are i.i.d.\ with the common law determined by 
\begin{equation}
\wt\BbbP\bigl(\wt T(e)>u\bigr)=  \wt\BbbP\otimes P\bigl(T^\OFF_0-U(1+T^\OFF_0)>u\bigr)
\end{equation}
where~$P$ is the law of  a uniform random variable~$U$~$[0,1]$ which (under $\wt\BbbP\otimes P$) is  independent of $T^\OFF_0$. The assumption  $T_0^\OFF\not\in L^{1/2}(\wt\BbbP)$  then forces $\wt\E (\wt T(e)^{1/2})=\infty$ and so, by the standard facts about sequences of i.i.d.\ random variables,
\begin{equation}
\label{E:4.9}
\forall\epsilon>0\colon\quad \frac1n\max_{0\le x\le \epsilon\sqrt n}\wt T(x,x+1)\,\underset{n\to\infty}\longrightarrow\,\infty\qquad\wt\BbbP\text{-a.s.}
\end{equation}
Since the random walk~$X$ cannot cross edge~$e$ before time~$\wt T(e)$, on the event that the maximum in \eqref{E:4.9} is larger than~$n$ we have $\max_{t\in[0,n]}X_t\le\epsilon\sqrt n$ almost surely. By symmetry, $X_t=o(\sqrt t)$ a.s.\ as~$t\to\infty$ thus showing that~$X$ is subdiffusive.
\end{proofsect}

\vbox{
\section*{Acknowledgments}
\nopagebreak\noindent
This work has been partially supported by NSF award DMS-1954343. 
}

\end{document}